\documentclass[11pt,amsfonts]{article}
\usepackage{graphicx,amssymb,eucal,mathrsfs}
\usepackage{latexsym,amsmath,epsfig,epic,eepic}
\usepackage[all]{xy}
\usepackage[english, french]{babel}

\newtheorem{theorem}{{\bf Th\'eor\`eme}}[section]
\newtheorem{main}{{\bf  Th\'eor\`eme}}
\newtheorem{proposition}[theorem]{{\bf Proposition}}

\newcommand{\RR}{\mathbb{\bf R}^{n}}

\newcommand{\R}{\mathbb{\bf R}}
\newcommand{\Z}{\mathbb{\bf Z}}

\newcommand{\G}{\mathscr{G}}

\newcommand{\tildem}{\widetilde{M}}
\newcommand{\identity}{1\hskip-2.5pt{\rm l}}

\title{Quelques plats pour la m\'etrique de Hofer}
\author{Pierre Py}

\begin{document}
\maketitle

\selectlanguage{english}

\begin{abstract}
We show, by an elementary and explicit construction, that the group of Hamiltonian diffeomorphisms of certain symplectic manifolds, endowed with Hofer's metric, contains subgroups quasi-isometric to Euclidean spaces of arbitrary dimension. (1991 MSC: 20F65, 37J05, 53D05)
\end{abstract}

\selectlanguage{french}

\section{Introduction}

Le but de ce texte est de prouver que le groupe des \textit{diff\'eomorphismes hamiltoniens} de certaines vari\'et\'es symplectiques, muni de la \textit{distance de Hofer}, contient des sous-groupes quasi-isom\'etriques \`a des espaces euclidiens de dimension arbitraire. Avant d'aller plus loin, rappelons quelques d\'efinitions. 

Consid\'erons une vari\'et\'e symplectique $(M,\omega)$ (connexe, \'eventuellement non compacte). Soit $H : M\times [0,1]\to \R$ une fonction lisse \`a support compact. Nous noterons souvent $H_{t}(x)=H(x,t)$. Le champ de vecteurs d\'ependant du temps $X_{H_{t}}$ d\'efini par la relation $\iota_{X_{H_{t}}}\omega=-dH_{t}$ s'int\`egre pour donner naissance \`a une isotopie $(f_{t})$ issue de l'identit\'e. Le diff\'eomorphisme $f_{1} : M \to M$ est le \textit{diff\'eomorphisme hamiltonien} engendr\'e par la fonction $H$. L'ensemble des diff\'eomorphismes ainsi obtenus forme un groupe, que l'on notera $\G$, qui est contenu dans le groupe des diff\'eomorphismes symplectiques de $M$. 

En 1990, Hofer \cite{hofer1} a d\'ecouvert que l'on pouvait munir le groupe $\G$ d'une remarquable distance, not\'ee $\rho$, qui est biinvariante, c'est-\`a-dire, invariante \`a la fois par les translations \`a droite et \`a gauche de $\G$. Rappelons-en la d\'efinition. L'oscillation d'une fonction \mbox{$F : M\to \R$}, not\'ee ${\rm osc}(F)$, est la quantit\'e
$${\rm max}(F)-{\rm min}(F).$$ 
Si $(H_{t})_{0\le t\le 1}$ est un hamiltonien d\'ependant du temps, nous pouvons d\'efinir la longueur de l'isotopie hamiltonienne $\{f_{t}\}$ engendr\'ee par $(H_{t})$ par :
$$\ell (\{f_{t}\})=\int_{0}^{1}{\rm osc}(H_{t})dt.$$ L'\textit{\'energie} d'un \'el\'ement $f\in \G$ est la quantit\'e 
$$\vert \vert f \vert \vert ={\rm inf}\, \ell (\{f_{t}\})$$ o\`u l'infimum porte sur toutes les isotopies hamiltoniennes dont le temps $1$ est $f$. On d\'efinit alors $\rho(f,g)=\vert \vert fg^{-1}\vert \vert$. Notons que le point crucial pour s'assurer que $\rho$ est une distance est d'\'etablir qu'un diff\'eomorphisme $f\in \G \setminus \{\identity \}$ a une \'energie strictement positive. Cela a \'et\'e prouv\'e par Hofer \cite{hofer1} dans le cas de $\R^{2n}$, par Polterovich \cite{pol1} pour les vari\'et\'es rationnelles et en toute g\'en\'eralit\'e par Lalonde et McDuff \cite{lmc} (voir aussi \cite{chekanov,oh1,schwarz,viterbo} pour d'autres preuves, dans diff\'erents cas particuliers). Nous reviendrons sur ce point au paragraphe suivant. Nous renvoyons le lecteur aux livres \cite{hoferzehnder,pol3} pour une introduction plus d\'etaill\'ee \`a ce sujet.

Bialy et Polterovich \cite{bialypolterovich} ont prouv\'e le r\'esultat suivant. Lorsque $M=\R^{2n}$, muni de sa structure symplectique standard, il existe un voisinage $U$ de l'identit\'e dans $\G$ (pour la topologie $C^{1}$) et un voisinage $V$ de l'origine dans l'espace vectoriel des fonctions $C^{\infty}$ \`a support compact sur $\RR$ (pour la topologie $C^{2}$) tels que $(U,\rho)$ et $(V,d_{{\rm osc}})$ soient isom\'etriques (o\`u $d_{{osc}}(F,G)={\rm osc}(F-G)$). On peut interpr\'eter ce fait en disant que $\G$ est \textit{localement plat}. Ce r\'esultat a depuis \'et\'e g\'en\'eralis\'e \`a d'autres vari\'et\'es \cite{lmc2,oh2}.

Par ailleurs, on sait maintenant pour une large classe de vari\'et\'es symplectiques que l'espace m\'etrique $(\G,\rho)$ est de diam\`etre infini (voir \cite{pol3} pour quelques r\'esultats et r\'ef\'erences). A l'oppos\'e du r\'esultat de Bialy et Polterovich pr\'ec\'edemment cit\'e, on peut donc s'int\'eresser, lorsque le diam\`etre de $\G$ est infini, \`a la g\'eom\'etrie \textit{\`a grande \'echelle} de $(\G,\rho)$. Dans cet esprit nous montrons le :

\begin{main}Supposons qu'il existe une sous-vari\'et\'e lagrangienne ferm\'ee $L$ plong\'ee dans $M$, v\'erifiant les deux conditions suivantes :

\noindent $\bullet$ l'application induite $\pi_{1}(L)\to \pi_{1}(M)$ entre les groupes fondamentaux de $L$ et $M$ est injective,

\noindent $\bullet$ il existe sur $L$ une m\'etrique riemannienne \`a courbure n\'egative ou nulle.

\noindent Alors, pour tout entier naturel $N$, il existe un morphisme $\phi : \Z^{N}\to \G$ ayant la propri\'et\'e suivante. Si $\vert \cdot \vert$est une norme fix\'ee sur $\R^{N}$, il existe une  constante strictement positive $C_{N}$ telle que :  
$$C_{N}^{-1}\vert x-y\vert \le \rho(\phi(x),\phi(y))\le C_{N}\vert x-y \vert, $$ pour tous $x,y$ de $\Z^{N}$.

\end{main}

Signalons \`a l'attention du lecteur anglophone, que les expressions ``n\'egative ou nulle" ou encore ``n\'egative", se traduisent par ``nonpositive". Citons quelques exemples de vari\'et\'es symplectiques  v\'erifiant les hypoth\`eses du th\'eor\`eme : 

\noindent $\bullet$ Le fibr\'e cotangent $T^{*}L$ (muni de sa structure symplectique canonique : $\omega=d(pdq)$), d'une vari\'et\'e ferm\'ee $L$ poss\'edant une m\'etrique riemannienne \`a courbure n\'egative ou nulle.  

\noindent $\bullet$ Une surface compacte orientable de genre strictement positif, munie d'une forme d'aire ; un produit de surfaces de genres strictement positifs.  

\noindent $\bullet$ Soit $V^{3}$ une vari\'et\'e ferm\'ee de dimension $3$ qui fibre sur le cercle : $\pi : V^{3}\to {\bf S}^{1}$. Supposons le genre de la fibre strictement positif. Notons $\theta_{1}$ la coordonn\'ee sur le cercle. Soit $\Omega$ une $2$-forme ferm\'ee sur $V^{3}$ qui soit non-d\'eg\'en\'er\'ee sur chaque fibre de $\pi$. Consid\'erons la vari\'et\'e $M=V^{3}\times {\bf S}^{1}$. Si $\theta_{2}$ d\'esigne la coordonn\'ee sur le second facteur de $M$, la forme $$\omega=\Omega+\pi^{*}(d\theta_{1})\wedge d\theta_{2}$$ est une forme symplectique sur $M$. Pour toute courbe ferm\'ee simple essentielle $\gamma$ contenue dans une fibre de $\pi$, nous obtenons un tore lagrangien incompressible $\gamma \times {\bf S}^{1}$ dans $M$.

{\bf Remarques.} 

\noindent $\bullet$ Des plongements quasi-isom\'etriques de groupes de type fini dans le groupe des diff\'eomor\-phismes hamiltoniens du disque $\mathbb{\bf D}^{2}=\{(x,y)\in \R^{2}, \vert x \vert ^{2}+\vert y \vert ^{2} <1\}$, muni de sa m\'etrique ${\rm L}^{2}$ (voir \cite{ak} \`a ce sujet), qui est invariante \`a gauche seulement, ont \'et\'e construits par Benaim et Gambaudo \cite{benaimgambaudo}, puis Crisp et Wiest \cite{crispwiest}. 

\noindent $\bullet$ Dans le cas o\`u $M$ est une surface ferm\'ee de genre strictement positif, les r\'esultats de \cite{pol3} permettent de plonger quasi-isom\'etriquement un espace de dimension infinie dans $\G$. Consid\'erons par exemple le tore $\R^{2}/\Z^{2}$, muni de la forme d'aire $dx\wedge dy$. Soit $E$ l'espace des fonctions de moyenne nulle sur $\R^{2}/\Z^{2}$, ne d\'ependant que de la coordonn\'ee $x\in \R/\Z$, muni de la norme $\vert H \vert_{\infty} = {\rm sup}_{x}\vert H(x)\vert$. Pour $H\in E$, notons $\psi(H)$ le temps $1$ du flot hamiltonien de $H$. Le th\'eor\`eme $7.2.C$ de \cite{pol3} assure que $\rho(\psi(H),\identity)\ge \vert H\vert_{\infty}$. Par ailleurs $\rho (\psi(H),\identity)\le {\rm osc}(H)\le 2\vert H\vert_{\infty}$ ; l'application $\psi : E \to \G$ est donc un plongement quasi-isom\'etrique.

\noindent $\bullet$ Pour construire le morphisme $\phi$, nous allons utiliser une id\'ee de Lalonde et \hbox{Polterovich} \cite{lalondepolterovich}. D'apr\`es un th\'eor\`eme classique de Weinstein, un voisinage $U$ de $L$ dans $M$ est symplectomorphe \`a un voisinage de la section nulle dans le fibr\'e cotangent $T^{*}L$. On peut supposer que $U$ est un fibr\'e en boules au-dessus de $L$. Nous pouvons donc consid\'erer des flots hamiltoniens sur $M$ qui, dans $U$ (ou une partie de $U$), co\"{i}ncident avec le flot g\'eod\'esique sur $L$ (pour une m\'etrique \`a courbure n\'egative fix\'ee). Notons $\widetilde{u}^{t} : T^{*}\widetilde{L}\to  T^{*}\widetilde{L}$ le flot g\'eod\'esique sur le rev\^etement universel $\widetilde{L}$ de $L$. Nous allons tirer parti du fait suivant, d\^u \`a la courbure n\'egative : si $K$ est un compact de $T^{*}\widetilde{L}$ qui \'evite la section nulle, on a $\widetilde{u}^{t}(K)\cap K= \emptyset$ pour $t$ suffisament grand.

\noindent $\bullet$ Nous verrons que, lorsque l'on choisit comme norme sur $\R^{N}$ la norme 
$$\vert (x_{1},\ldots , x_{N}) \vert = \sum_{k=1}^{N}\vert x_{k}\vert $$ la constante $C_{N}$ que nous obtenons converge exponentiellement vite vers l'infini lorsque $N$ tend vers l'infini. 

\noindent $\bullet$ Il serait int\'eressant de savoir si un r\'esultat analogue est vrai lorsque $M$ est la sph\`ere $\mathbb{\bf S}^{2}$. Dans ce cas, et \`a notre connaissance, les seules mani\`eres d'obtenir des bornes inf\'erieures arbitrairement grandes sur la distance de Hofer proviennent de \cite{entpol,pol2}. Leonid Polterovich m'a indiqu\'e que, dans le cas o\`u $M$ est le disque $\mathbb{\bf D}^{2}$, les r\'esultats de \cite{entpol} permettent de prouver un r\'esultat analogue au th\'eor\`eme $1$.

Dans la seconde partie de ce texte, nous rappelons quelques faits classiques de topologie symplectique, puis nous prouvons le th\'eor\`eme dans la troisi\`eme partie. 

\section{In\'egalit\'e entre \'energie et capacit\'e}

Si $A$ est une partie de $M$, rappelons que la \textit{capacit\'e de Gromov} de $A$, not\'ee $c(A)$, est la quantit\'e : 
$${\rm sup}\{\pi r^{2}, {\rm il \, existe\, un\, plongement\,  symplectique\, } B^{2n}(r)\to {\rm Int}(A)\}.$$
Ici, $B^{2n}(r)$ d\'esigne la boule euclidienne de rayon $r$ dans $\R^{2n}$ muni de sa structure symplectique standard (et ${\rm dim} M= 2n$). Cette notion a \'et\'e introduite dans \cite{gromov}. Lalonde et McDuff \cite{lmc} ont \'etabli le r\'esultat suivant. Si $f$ est un \'el\'ement de $\G$ et $A$ une partie de $M$ qui est disjointe d'elle-m\^eme par $f$, c'est-\`a-dire qui v\'erifie $f(A)\cap A=\emptyset$, alors l'\'energie de $f$ est minor\'ee par la moiti\'e de la capacit\'e de $A$ : 
$$\vert \vert f \vert \vert \ge \frac{1}{2}c(A).$$ Rappelons que c'est cette in\'egalit\'e qui permet d'\'etablir que la distance de Hofer est non-d\'eg\'en\'er\'ee : si $f$ est un diff\'eomorphisme diff\'erent de l'identit\'e, on peut trouver un ouvert de $M$ qui est disjoint de lui-m\^eme par $f$. Puisque tout ouvert non-vide a une capacit\'e strictement positive, l'\'energie de $f$ est non-nulle. Une telle in\'egalit\'e avait \'et\'e prouv\'ee (sans le facteur $\frac{1}{2}$) par Hofer dans $\R^{2n}$ \cite{hofer2} (voir aussi \cite{fgs}).

Expliquons maintenant comment nous allons utiliser cette in\'egalit\'e. Notons d'abord que toutes les vari\'et\'es symplectiques qui v\'erifient les hypoth\`eses de notre th\'eor\`eme ont un groupe fondamental infini, et donc un rev\^etement universel $p : \widetilde{M}\to M$ non compact. Fixons un \'el\'ement $f$ de $\G$. Consid\'erons une isotopie hamiltonienne $f_{t} : M \to M$, engendr\'ee par un hamiltonien $(H_{t})$, telle que $f_{1}=f$. Notons $\widetilde{f}_{t} : \widetilde{M} \to \widetilde{M}$ le relev\'e de l'isotopie $(f_{t})$ issu de l'identit\'e. C'est une isotopie hamiltonienne engendr\'ee par la fonction $H_{t}\circ p$ (\`a support non compact). Nous obtenons ainsi un relev\'e  $\widetilde{f}_{1} : \widetilde{M}\to \widetilde{M}$ de $f$. Suivant \cite{lalondepolterovich}, nous appellerons \textit{relev\'e admissible} de $f$ tout relev\'e ainsi obtenu. 

\begin{proposition} \cite{lalondepolterovich} Soit $c>0$ et $f\in \G$. Supposons que pour tout relev\'e admissible $\widetilde{f}$ de $f$, il existe une partie $A$ de $\widetilde{M}$ de capacit\'e sup\'erieure ou \'egale \`a $c$ telle que $\widetilde{f}(A)\cap A = \emptyset$. Alors $\vert \vert f\vert \vert \ge \frac{c}{2}$ 
\end{proposition}

\noindent {\bf Remarque :} si $M$ est non compacte, puisque nous ne consid\'erons que des isotopies hamiltoniennes sur $M$ \`a support compact, il est clair qu'il existe un unique relev\'e admissible. Si $M$ est compacte, une cons\'equence facile de la (difficile) conjecture d'Arnold (voir \cite{fukayaono,liutian}) est que l'application d'\'evaluation
$$\pi_{1}(\G,\identity)\to \pi_{1}(M,x_{0})$$ a une image triviale. Ceci implique que tout diff\'eomorphisme hamiltonien de $M$ poss\`ede un {\bf unique} relev\'e admissible. Ainsi, pour appliquer la proposition ci-dessus, il suffit de v\'erifier l'hypoth\`ese pour un seul diff\'eomorphisme de $\widetilde{M}$. Cependant, pour garder \`a cette article un caract\`ere \'el\'ementaire, nous n'utiliserons pas ce fait. 

\noindent \textit{Preuve de la proposition} : consid\'erons une isotopie hamiltonienne $(f_{t})$ sur $M$, engendr\'ee par la fonction \`a support compact $H_{t}$, telle que $f_{1}=f$. Soit $A$ un compact contenant la r\'eunion des supports des fonctions $H_{t}$. Si $M$ est compacte, nous supposerons $H_{t}$ normalis\'ee par la condition $\int_{M}H_{t} \, \omega^{n}=0$, pour tout $t$. Soit $(\widetilde{f}_{t})$ l'isotopie de $\tildem$ engendr\'ee par la fonction $H_{t}\circ p$. Par hypoth\`ese on peut trouver un compact $K$ de $\tildem$ de capacit\'e sup\'erieure \`a $c-\epsilon$, qui est disjoint de lui-m\^eme par $\widetilde{f}_{1}$. Soit $B$ une boule de $\tildem$ qui contient $\cup_{t\in [0,1]}f_{t}(K)$ et telle que la projection $p : B \to A$ soit surjective ; et $\varphi : \tildem \to [0,1]$ une fonction \`a support compact valant $1$ sur $B$. L'isotopie hamiltonienne engendr\'ee par la fonction \`a support compact d\'efinie par $G_{t}(x)=\varphi(x)H_{t}(p(x))$ disjoint $K$ de lui-m\^eme. Nous obtenons donc, d'apr\`es l'in\'egalit\'e entre \'energie et capacit\'e, $\int_{0}^{1}{\rm osc}(G_{t})dt\ge \frac{c-\epsilon}{2}$. Puisque  ${\rm osc}(G_{t})={\rm osc}(H_{t})$, on obtient l'estimation voulue. \hfill $\Box$

Une cons\'equence classique de la preuve ci-dessus est qu'un relev\'e admissible d'un diff\'eomorphisme hamiltonien de $M$ ne peut disjoindre d'elle-m\^eme une partie de $\tildem$ de capacit\'e infinie.

\section{Preuve du th\'eor\`eme}

Fixons une m\'etrique riemannienne \`a courbure n\'egative ou nulle $g$ sur $L$. Nous noterons $d_{g}$ la distance induite par $g$ sur le rev\^etement universel $\widetilde{L}$ de $L$. Quitte \`a multiplier la m\'etrique $g$ par une constante, on peut supposer qu'il existe un voisinage $U$ de $L$ dans $M$ et un diff\'eomorphisme symplectique 
$$\theta : T^{*}L(\sqrt{3})\to U,$$ o\`u $T^{*}L(\sqrt{3})=\{(q,p)\in T^{*}L, \vert p\vert_{q}^{2}< 3\}$. Fixons d\'esormais un entier naturel $N$. Notons $A_{i}$ ($1\le i \le 2^{N}$) la partie suivante de $T^{*}L(\sqrt{3})$ :

$$\{(q,p), 1+\frac{i-1}{2^{N}}\le \vert p\vert_{q}^{2}\le 1+\frac{i-1}{2^{N}}+\frac{1}{2^{N+1}}\}.$$ Nous avons repr\'esent\'e en noir sur la figure 1, dans le cas o\`u $N=2$, la trace des ensembles $A_{1},\ldots , A_{4}$ sur une fibre de la projection $T^{*}L\to L$.

Il sera plus commode d'indexer en fait les ensembles $(A_{i})_{1\le i \le 2^{N}}$ par $\{\pm 1\}^{N}$. Pour cela nous fixons une bijection entre $\{\pm 1\}^{N}$ et $\{1,\dots , 2^{N}\}$ : 

$$I=(I_{1},\ldots , I_{N})\in \{\pm 1\}^{N}\mapsto i(I)\in \{1,\ldots , 2^{N}\}.$$
Soit, pour $1\le k \le N$, $\varphi_{k} : [0,3]\to \R$ une fonction de classe $C^{\infty}$, ayant les propri\'et\'es suivantes : 

\noindent $\bullet$ l'application $\varphi_{k}$ est nulle en dehors de l'intervalle $[\frac{1}{2}, \frac{5}{2}]$. 

\noindent $\bullet$ Si $I_{k}=1$ et $s\in [1+\frac{i(I)-1}{2^{N}},1+\frac{i(I)-1}{2^{N}}+\frac{1}{2^{N+1}}]$, $\varphi_{k}(s)=s$ ; si $I_{k}=-1$ et $s\in [1+\frac{i(I)-1}{2^{N}},1+\frac{i(I)-1}{2^{N}}+\frac{1}{2^{N+1}}]$, $\varphi_{k}(s)=-s$.

\noindent Enfin si $(q,p)\in T^{*}L(\sqrt{3})$, on pose $H_{k}(q,p)=\frac{1}{2}\varphi_{k}(\vert p\vert_{q}^{2})$. Gr\^ace au diff\'eomorphisme $\theta$ on peut voir $H_{k}$ comme une fonction sur $U$, que l'on prolonge par $0$ en dehors de $U$ pour obtenir une fonction lisse sur $M$. Les flots hamiltoniens $\phi_{H_{k}}^{t}$ associ\'es aux fonctions $H_{1},\ldots, H_{N}$ commutent et d\'efinissent une action de $\R^{N}$ sur $M$. Nous d\'efinissons un morphisme $\phi : \Z^{N} \to \G$ par : 
$$\phi(a=(a_{1}, \ldots , a_{N}))=\prod_{k=1}^{N} \phi_{H_{k}}^{a_{k}}.$$

\noindent Nous avons bien s\^ur : 
$$\rho(\phi(a),\phi(b))=\vert \vert \prod_{k=1}^{N}\phi_{H_{k}}^{a_{k}-b_{k}}\vert \vert \le C\cdot (\sum_{k=1}^{N}\vert a_{k}-b_{k}\vert),$$ o\`u l'on a not\' e $C={\rm max}_{1\le k \le N}\vert \vert \phi_{H_{k}}^{1}\vert \vert $. Pour prouver le th\'eor\`eme, nous devons \'etablir une minoration de la forme $$\vert \vert \phi(a) \vert \vert \ge \epsilon_{N}\cdot (\sum_{k=1}^{N}\vert a_{k}\vert)$$
pour tout $a \in \Z^{N}$, pour une certaine constante $\epsilon_{N}$. On peut ensuite prendre $C_{N}={\rm max}(\epsilon_{N}^{-1},C)$. 

\begin{figure}
\begin{center}

\setlength{\unitlength}{0.00087489in}
\begingroup\makeatletter\ifx\SetFigFont\undefined%
\gdef\SetFigFont#1#2#3#4#5{%
  \reset@font\fontsize{#1}{#2pt}%
  \fontfamily{#3}\fontseries{#4}\fontshape{#5}%
  \selectfont}%
\fi\endgroup%
{\renewcommand{\dashlinestretch}{30}
\begin{picture}(1318,1333)(0,-10)
\put(659,659){\whiten\ellipse{866}{866}}
\put(659,659){\ellipse{866}{866}}
\put(659,659){\blacken\ellipse{816}{816}}
\put(659,659){\ellipse{816}{816}}
\put(659,659){\whiten\ellipse{764}{764}}
\put(659,659){\ellipse{764}{764}}
\put(659,659){\blacken\ellipse{708}{708}}
\put(659,659){\ellipse{708}{708}}
\put(659,659){\whiten\ellipse{654}{654}}
\put(659,659){\ellipse{654}{654}}
\put(659,659){\blacken\ellipse{594}{594}}
\put(659,659){\ellipse{594}{594}}
\put(659,659){\whiten\ellipse{540}{540}}
\put(659,659){\ellipse{540}{540}}
\put(659,659){\blacken\ellipse{490}{490}}
\put(659,659){\ellipse{490}{490}}
\put(659,659){\whiten\ellipse{438}{438}}
\put(659,659){\ellipse{438}{438}}
\put(659,659){\ellipse{1302}{1302}}
\end{picture}
}
\caption{}

\end{center}
\end{figure}

Fixons donc $a=(a_{1}, \ldots, a_{N})\in \Z^{N}- \{0\}$. On peut choisir $I\in \{\pm 1\}^{N}$ tel que $I_{k}a_{k}\ge 0$ pour tout $k$ de $\{1,\ldots , N\}$. Si $(q,p)\in A_{i(I)}$, alors $$(\sum_{k=1}^{N}a_{k}H_{k})(q,p)=\frac{l}{2}\vert p\vert_{q}^{2}$$ (o\`u $l:=\sum_{k=1}^{N}\vert a_{k}\vert $). Notons $u^{t} : T^{*}L \to T^{*}L$ le flot g\'eod\'esique (pour la m\'etrique $g$) : c'est le flot hamiltonien associ\'e \`a la fonction $E(q,p)=\frac{1}{2}\vert p\vert^{2}_{q}$. Le flot hamiltonien engendr\'e par la fonction $\sum_{k=1}^{N}a_{k}H_{k}$ co\"{i}ncide donc avec le flot $(u^{lt})$ sur $A_{i(I)}$.  

Choisissons une composante connexe $\widetilde{U}$ de l'image inverse de $U$ dans $\widetilde{M}$. Puisque $L$ est incompressible dans $M$, $\widetilde{U}$ est symplectiquement diff\'eomorphe \`a $T^{*}\widetilde{L}(\sqrt{3})$. Nous choisissons \'egalement un point base $q_{0}\in \widetilde{L}$, et notons $B(q_{0},R)$ la boule ouverte de rayon $R$ centr\'ee en $q_{0}$.

Notons $\widetilde{\phi}(a) : \widetilde{M}\to \tildem$ le relev\'e admissible de $\phi(a)$ d\'etermin\'e par l'isotopie engendr\'ee par la fonction $\sum_{k=1}^{N}a_{k}H_{k}$. Un autre relev\'e admissible (hypoth\'etique !!! d'apr\`es la remarque faite plus haut) serait de la forme $T\circ \widetilde{\phi}(a)$ o\`u $T$ est un \'el\'ement du groupe fondamental de $M$. Notons que $T$ appartient n\'ecessairement au groupe fondamental de $L$, sinon $T\circ \widetilde{\phi}(a)$ disjoindrait $\widetilde{U}$ de lui-m\^eme. C'est impossible car $\widetilde{U}$ est de capacit\'e infinie (ceci se d\'eduit, par exemple, de la proposition 3.1).

Posons $R=\frac{l}{4}$. Soit $\widetilde{A}_{i,R}$ la partie suivante de $T^{*}\widetilde{L}(\sqrt{3})$ :
$$\{(q,p), q\in B(q_{0},R),  1+\frac{i-1}{2^{N}}\le \vert p\vert_{q}^{2}\le 1+\frac{i-1}{2^{N}}+\frac{1}{2^{N+1}}\}.$$ Si $(q,p)\in \widetilde{A}_{i,R}$, \'ecrivant $\widetilde{\phi}(a)(q,p)=(q',p')$,  nous avons $d_{g}(q,q')\ge l$. Ceci assure que $\widetilde{\phi}(a)(\widetilde{A}_{i,R})\cap \widetilde{A}_{i,R}=\emptyset$.

\noindent Consid\'erons maintenant un autre relev\'e admissible de la forme $T\circ \widetilde{\phi}(a)$ ($T\in \pi_{1}(L)- \{1\}$). Supposons que $T\circ \widetilde{\phi}(a)(\widetilde{A}_{i,R})$ rencontre $\widetilde{A}_{i,R}$. Il existe alors $(q,p)\in \widetilde{A}_{i,R}$ tel que $\widetilde{\phi}(a)(q,p)=(q',p')\in T^{-1}(\widetilde{A}_{i,R})$. On obtient alors $$d_{g}(T(q_{0}),q_{0})\ge d_{g}(T(q_{0}),T(q'))-R\ge d_{g}(q',q)-2R\ge \frac{l}{2}.$$ Si $\nu$ et $C$ sont des constantes strictement positives, nous noterons $$\Lambda_{\nu,C}=\{(q,p),q\in B(q_{0},C), \vert p\vert_{q}<\nu\}.$$ Rappelons que l'on a $H_{k}(q,p)=0$ (pour tout $k$), d\`es que $\vert p \vert_{q}^{2}<\frac{1}{2}$. Alors, si $(q,p)\in \Lambda_{\frac{1}{\sqrt{2}},\frac{R}{2}}$, $T\circ \widetilde{\phi}(a)(q,p)=T(q,p)=(q',p')$ v\'erifie : 
$$d_{g}(q',q_{0})\ge d_{g}(T(q_{0}),q_{0})-d_{g}(q,q_{0})\ge \frac{l}{2}-\frac{l}{8}=\frac{3l}{8}.$$ Donc $T\circ \widetilde{\phi}(a)(\Lambda_{\frac{1}{\sqrt{2}},\frac{R}{2}})\cap(\Lambda_{\frac{1}{\sqrt{2}}, \frac{R}{2}})=\emptyset$.

En r\'esum\'e, tout relev\'e admissible de $\phi(a)$ disjoint d'elle-m\^eme une partie de $\tildem$ de capacit\'e sup\'erieure ou \'egale \`a ${\rm min}(c(\widetilde{A}_{i,R}),c(\Lambda_{\frac{1}{\sqrt{2}}, \frac{R}{2}}))$. D'apr\`es la proposition 2.1, nous avons $\vert \vert \phi(a)\vert \vert \ge \frac{1}{2}{\rm min}(c(\widetilde{A}_{i,R}),c(\Lambda_{\frac{1}{\sqrt{2}}, \frac{R}{2}}))$. Pour conclure la preuve du th\'eor\`eme, il nous reste \`a obtenir une minoration, lin\'eaire en $R$, des capacit\'es de $\widetilde{A}_{i,R}$ et $\Lambda_{\frac{1}{\sqrt{2}},\frac{R}{2}}$. La preuve de la proposition suivante m'a \'et\'e sugg\'er\'ee par Jean-Claude Sikorav. 

\begin{proposition} Il existe une constante $\varepsilon >0$ telle que ${\rm min}(c(\widetilde{A}_{i,R}),c(\Lambda_{\frac{1}{\sqrt{2}},\frac{R}{2}}))\ge \varepsilon R$.

\end{proposition} 

\noindent \textit{Preuve} : on commence par ramener l'estimation de la capacit\'e de $\widetilde{A}_{i,R}$ \`a celle d'un ensemble de la forme $\Lambda_{\alpha,R}=\{(q,p), q\in B(q_{0},R), \vert p\vert_{q}<\alpha \}$.

\noindent Soit $V : \widetilde{L} \to \R$ une fonction de classe $C^{\infty}$ telle que $\vert \vert dV(q)\vert \vert =\sqrt{1+\frac{i-1}{2^{N}}+\frac{1}{2^{N+2}}}$ pour $q\in B(q_{0},2R)$. Pour cela fixons un point $q_{\infty}\in \widetilde{L}$ tel que $d_{g}(q_{0},q_{\infty})\ge 10^{7}\cdot R$,  et prenons pour $V$ la fonction $$\sqrt{1+\frac{i-1}{2^{N}}+\frac{1}{2^{N+2}}}d_{g}(\cdot ,q_{\infty}),$$ multipli\'ee par une fonction plateau qui s'annule au voisinage de $q_{\infty}$. Notons $T_{V}(q,p)=(q,p-dV(q))$ ($(q,p)\in T^{*}\widetilde{L}$). On v\'erifie ais\'ement que si $q\in B(q_{0},R)$ et $\vert p \vert_{q}< (10\cdot 2^{N+2})^{-1}$, alors $(q,dV(q)+p)\in \widetilde{A}_{i,R}$. L'application $T_{V}$ \'etant un diff\'eomorphisme symplectique de $T^{*}\widetilde{L}$, on a donc 
$$c(\widetilde{A}_{i,R})=c(T_{V}(\widetilde{A}_{i,R}))\ge c(\Lambda_{(10\cdot 2^{(N+2)})^{-1},R}).$$ Il nous reste maintenant \`a obtenir une minoration de la capacit\'e de $\Lambda_{\alpha, R}$. L'application $G : T_{q_{0}}\widetilde{L}\times T_{q_{0}}^{*}\widetilde{L}\to T^{*}\widetilde{L}$ d\'efinie par : $$(v\in T_{q_{0}}\widetilde{L}, \eta \in T_{q_{0}}^{*}\widetilde{L})\mapsto (exp_{q_{0}}(v),\eta \circ\left(D exp_{q_{0}}(v) \right)^{-1})$$ est un diff\'eomorphisme symplectique : c'est l'application induite entre les fibr\'es cotangents de $T_{q_{0}}\widetilde{L}$ et $\widetilde{L}$ par le diff\'eomorphisme $exp_{q_{0}} : T_{q_{0}}\widetilde{L}\to \widetilde{L}$. Puisque la m\'etrique $g$ est \`a courbure n\'egative ou nulle, l'application $(D exp_{q_{0}}(v))^{-1} : (T_{exp_{q_{0}}(v)}\widetilde{L},g_{exp_{q_{0}}(v)}) \to (T_{q_{0}}\widetilde{L},g_{q_{0}})$ est de norme major\'ee par $1$. Il en est de m\^eme pour sa transpos\'ee. On a donc 
$$\{(v,\eta), \vert v \vert _{q_{0}}< R, \vert \eta \vert_{q_{0}} < \alpha \} \subset G^{-1}(\Lambda_{\alpha , R}).$$
 L'application lin\'eaire symplectique 
 $$(q,p)\mapsto (\sqrt{\frac{R}{\alpha}}q,\sqrt{\frac{\alpha}{R}}p) $$ envoie la boule euclidienne $B(0,\sqrt{R\alpha})$ dans $\{(v,\eta), \vert v \vert _{q_{0}}< R, \vert \eta \vert_{q_{0}} < \alpha\}$. Nous obtenons donc bien l'in\'egalit\'e $c(\Lambda_{\alpha,R})\ge \pi R\alpha$. Finalement : $${\rm min} (c(\widetilde{A}_{i,R}),c(\Lambda_{\frac{1}{\sqrt{2}},\frac{R}{2}}))\ge {\rm min}(\frac{\pi R}{2\sqrt{2}},\frac{\pi R}{10\cdot 2^{N+2}})\ge \frac{\pi R}{10\cdot 2^{N+2}}\, .$$\hfill $\Box$

\medskip

\begin{flushleft}

\nopagebreak{ Pierre Py\\ ANR-06-BLAN-0030 \\ Unit\'e de Math\'ematiques Pures et Appliqu\'ees \\
UMR 5669 CNRS \\
\'Ecole Normale Sup\'erieure de Lyon\\ 46, All\'ee d'Italie\\ 69364 Lyon 
Cedex 07 \\ FRANCE\\ Pierre.Py@umpa.ens-lyon.fr}
\end{flushleft}

\end{document}